\documentclass[11pt,reprint]{revtex4-1}

\usepackage{graphicx,epsfig,amssymb,amsbsy,amsmath,amsthm}
\newcommand{\bb}[1]{\left({#1}\right)}					
\newcommand{\sq}[1]{\left[#1\right]}						
\newcommand{\cc}[1]{\left\{#1\right\}}					
\newcommand{\op}[1]{\mathcal{#1}}
\newcommand{\ord}[1]{{\sf O}\bb{#1}}					
\newcommand{\abs}[1]{\left|#1\right|}					

\newcommand{\sfrac}[2]{\mbox{$\frac{#1}{#2}$}}	
\newcommand{\hf}{\mbox{$\frac12$}}

\newcommand{\lambdab}{{\mbox{$\lambda$}}}

\newcommand{\sw}{v}

\newtheorem{example}{Example}

\newcommand{\sign}{\operatorname{sign}}
\newcommand{\step}{\operatorname{step}}
\newcommand{\erf}{\operatorname{Erf}}

\newcommand{\x}{X}

\newcommand{\f}{F}

\newcommand{\sa}{\alpha}
\renewcommand{\sb}{\beta}
\renewcommand{\sc}{\gamma}

\newcommand{\fref}[1]{figure~\ref{#1}}
\newcommand{\Fref}[1]{Figure~\ref{#1}}

\newcommand{\eref}[1]{(\ref{#1})}

\newcommand{\sref}[1]{section~\ref{#1}}
\newcommand{\Sref}[1]{Section~\ref{#1}}

\def\eps{\varepsilon}

\begin{document}

\title{Smoothing tautologies, hidden dynamics, and sigmoid asymptotics\\ for piecewise smooth systems}
\author{Mike R. Jeffrey}\address{Engineering Mathematics, University of Bristol, Merchant Venturer's Building, Bristol BS8 1UB, UK, email: mike.jeffrey@bristol.ac.uk}
\date{\today}

\begin{abstract}
Switches in real systems take many forms, such as impacts, electronic relays, mitosis, and the implementation of decisions or control strategies. To understand what is lost, and what can be retained, when we model a switch as an instantaneous event, requires a consideration of so-called {\it hidden} terms. These are asymptotically vanishing outside the switch, but can be encoded in the form of nonlinear switching terms. A general expression for the switch can be developed in the form of a series of sigmoid functions. We review the key steps in extending the Filippov's method of {\it sliding modes} to such systems. We show how even slight nonlinear effects can hugely alter the behaviour of an electronic control circuit, and lead to `hidden' attractors inside the switching surface. 
\end{abstract}

\maketitle



There are increasingly general and powerful methods to study the dynamical effects of switches in otherwise smoothly evolving dynamical systems, whether a phase transition, a trade or policy decision, an electrical relay, a mechanical actuator or a biochemical valve. This paper is an attempt to highlight and deal with certain ambiguities of non-uniqueness of {\it piecewise smooth} systems like these. Attempts are often made to show the robustness and generality of discontinuous models by smoothing or {\it regularizing} a discontinuity, a process that is somewhat tautological because, as we show, entirely different dynamics can be obtained with different choices of regularization. Fortunately this can be resolved using nonlinear switching terms, and we introduce a general method by developing sigmoid series expressions for discontinuous vector fields. We apply these ideas to a model of an elementary electronic control circuit, and give an example of how novel attractors may arise inside the switching surface. 

\section{Introduction}\label{sec:intro}

Many phenomena, in either our passive attempts to describe nature, or our active attempts to control it, are mixtures of steady behaviours and sharp transitions. 
The steady regimes are usually relatively easy to model. 
The transitions are often complicated and challenging to model, but occupy fleeting instants of time which make it tempting to simplify them, in the extreme, as discrete events in an otherwise smooth system. 
The effect of this, and how far we can develop the mathematics of a piecewise determined system, is explored here. 

Modeling a switch by means of piecewise smooth functions is simple and convenient, requiring no detailed in-depth knowledge of the physical laws during switching, but it risks destroying information about the transition process. The resulting simplicity has inspired an extensive theory of discontinuity-induced bifurcations and singularities. Much of the theoretical development is impressively powerful, and even highly rigorous (for reviews see \cite{f88,c02,physDspecial,machina2013siads,bc08}). Its generality, its applicability to {\it real} world switching phenomena, however, remains poorly understood (see e.g. \cite{guck-review}). 
We can, for example, derive piecewise smooth models as the limiting cases of smooth multiple timescale models, collapsing the fast switching dynamics into a single instant \cite{j15douglas}, or as describing the average dynamics of stochastic \cite{simpson12} or hysteretic \cite{f88,u92} processes. But the limits are typically singular, and qualitatively different systems have the same piecewise smooth limit, resulting in ambiguity about the systems they describe. 

A resolution to this has emerged recently. 
It seems that generality in piecewise smooth systems requires a notion of hidden dynamics that is observable only {\it during} a switch and not outside it, encoded in nonlinear terms that represent the remnants of asymptotically small effects that vanish in the limit of taking a piecewise smooth model. We shall present a few examples where nonlinear terms entirely alter the local and global behaviour of a system, including an application to electrical control systems. 

In engineering and biological applications particularly, an increasing number of examples of `nonlinear' switches are becoming known. 
These are systems whose dynamical equation contain nonlinear combinations of steep or discontinuous `sigmoid' functions, each sigmoid being an empirical representation of a particular switching action. In Boolean switching models of protein dynamics regulated via Hill functions, they are very common \cite{hill,abouthill08,abouthill12}. In \cite{j13error}, it was shown that nonlinear switching can be used to model static friction without needing to augment Coulomb friction with a velocity dependence. In general, nonlinear switching can be used to model effects such as overshoot or lag, and can resolve deterministic ambiguities that arise when multiple switches coincide. It seems reasonable to hope that a deeper understanding of them will permit greater use for modeling across the physical, social, and biological sciences. 

%
%
%
%

In \sref{sec:convex} we recount the standard method for handling a discontinuity in a dynamical system, and give examples of its ambiguity under smoothing. 
In \sref{sec:smooth} we discuss ambiguities in smoothing more generally. 
In \sref{sec:sigseries} we make these ambiguities explicit, and find how to express them in the piecewise smooth system. We apply these ideas to a circuit model in \sref{sec:circuit}. 
In \Sref{sec:non} we outline the method for solving a dynamical system at a discontinuity, followed by an example of the potential for novel dynamics, in the form of a hidden attractor, in \sref{sec:ueda}. 
A few closing remarks are made in \sref{sec:conc}.


\section{Nonsmooth dynamics -- the convex approach}\label{sec:convex}

Consider a piecewise smooth dynamical system
\begin{equation}\label{ns1}
\frac{dx}{dt}=f(x)=\left\{\begin{array}{lll}f_+(x)&\rm for&\sw(x)>0\;,\\f_-(x)&\rm for&\sw(x)<0\;.\end{array}\right.
\end{equation}
This paper is about how we extend such a model across $\sw(x)=0$, and the issues that arise in doing so.

In the standard Filippov approach \cite{f88} prominent in variable structure control \cite{u77,u92,es98} and in piecewise smooth dynamical systems theory \cite{f88,bc08}, when a system switches between two systems as in \eref{ns1}, we form their convex combination
\begin{equation}\label{ffil}
\frac{dx}{dt}\;=\;f(x;\lambda)\;=\;\frac{f_+(x)+ f_-(x) }2+\frac{f_+(x)- f_- (x)}2\lambda
\end{equation}
(or a similar form $\sfrac{dx}{dt}=\sq{f_+(x)+ f_-(x)}+u\sq{f_+(x)- f_- (x)}$ with $\lambda=2u-1$), where
\begin{equation}\label{sign}
\lambda\in\left\{\begin{array}{lll}\sign\bb{\sw}&\rm if&\sw\neq0\;,\\\sq{-1,+1}&\rm if&\sw=0\;,\end{array}\right.
\end{equation}
thus
\begin{equation}
f(x;+1)\equiv f_+(x)\qquad\mbox{and}\qquad f(x;-1)\equiv f_-(x)\;.
\end{equation}
The standard approach then seeks so-called {\it sliding modes} which satisfy $\sfrac{d\;}{dt}\sw=0$ on $\sw=0$. Taking coordinates $x=(x_1,x_2,...,x_n)$ such that $x_1=\sw$, and writing $f=(f_1,f_2,...,f_n)$, the dynamics of sliding modes is therefore given by 
\begin{equation}\label{blowslidingf}
\left.\begin{array}{rcl}0&=&f_1({x};\lambda^s)\\
\frac{d\;}{dt}\bb{x_2,...,x_n}&=&\bb{f_2({x};\lambdab^s),...,f_n({x};\lambdab^s)}\end{array}\right\}\quad{\rm on}\;\;x_1=0\;
\end{equation}
for some $\lambda^s\in[-1,+1]$. If no solutions of the sliding mode problem exist for $\lambda^s\in[-1,+1]$, then the flow of \eref{ffil} crosses through $\sw=0$ transversally. 

This definition of sliding and crossing at a switching surface is sometimes justified by smoothing (or `regularizing') the combination \eref{ffil}, by replacing $\lambda$ with a smooth sigmoid function, and showing that slow dynamics of the smoothed system contains invariant manifold dynamics equivalent to the sliding mode dynamics \eref{blowslidingf}, (for theorems regarding this equivalence, which remains important in spite of what follows, see \cite{st96,ts11,j15douglas}). 
Unfortunately such an argument is tautologous. Smoothing simply preserves the dynamics we have imposed already by writing the expression \eref{ffil}, and we shall see that by changing \eref{ffil} we can not only obtain different dynamics, but that the difference persists under smoothing. Two examples illustrate the ambiguities of smoothing a discontinuity very clearly, and while almost trivially simple, they help motivate the general arguments that follow. 

\begin{figure}[h!]\centering\includegraphics[width=0.45\textwidth]{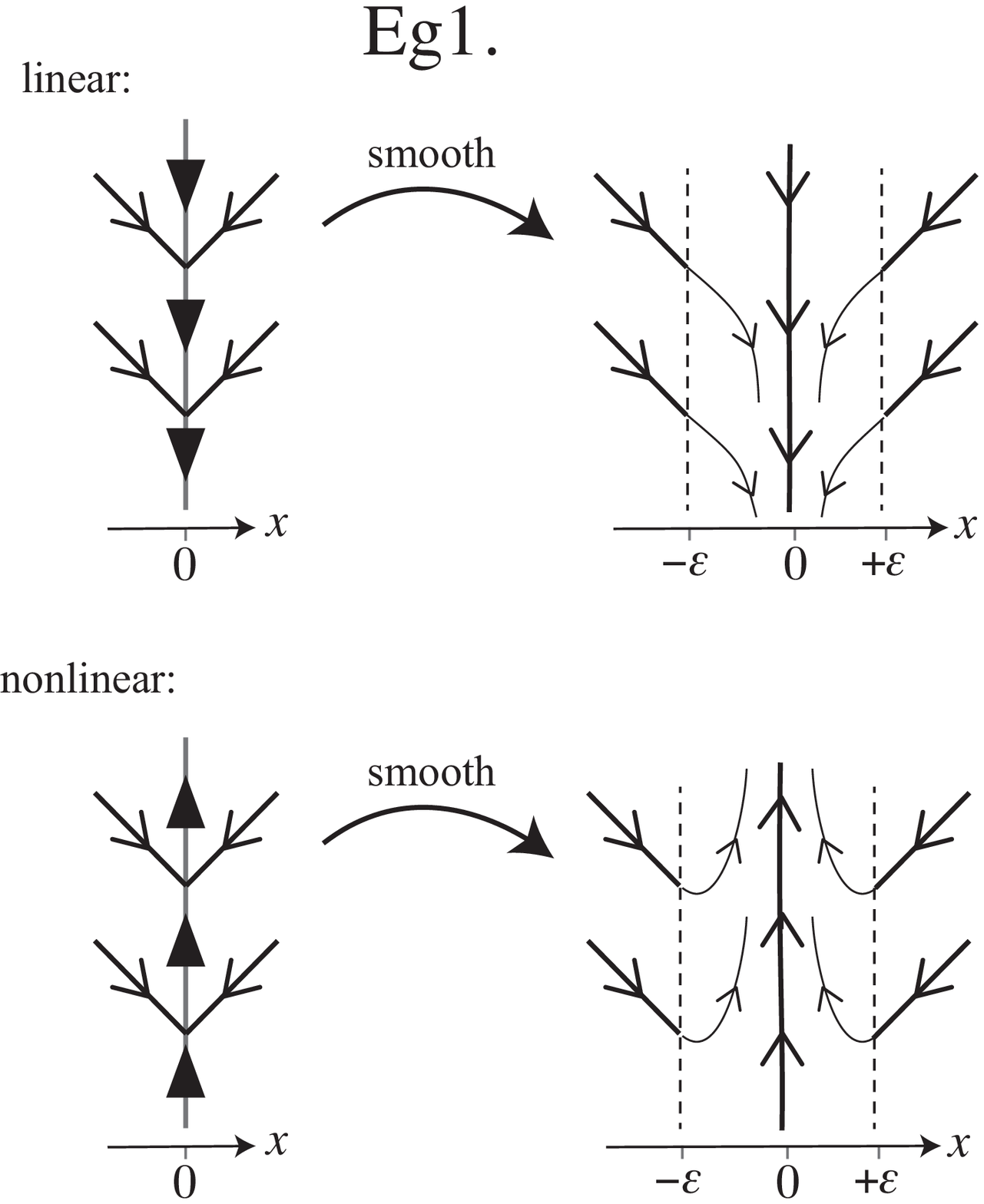}
\vspace{-0.3cm}\caption{\sf The system $\sfrac{d\;}{dt}(x_1,x_2)=({\rm sign}(x_1),1-2({\rm sign}(x_1))^2)$ viewed ignoring the nonlinear dependence on $\sign(x_1)$ (left) or respecting it (right), in each case showing the discontinuous system and its smoothing. The phase portraits as shown correspond to equations \eref{eg1f}, \eref{eg1sfs}, \eref{snonfil}, \eref{snon1}, respectively. 
}\label{fig:eg1}\end{figure}

\begin{example}
Let us start with a smooth system that we wish to study, 
\begin{equation}
\sfrac{d\;}{dt}(x_1,x_2)=(\varphi_\eps(x_1),1-2\varphi_\eps^2(x_1))\;,\label{snon1}
\end{equation} 
where $\eps$ is a small positive parameter and $\varphi_\eps(x_1)$ is a sigmoid function such that $\varphi_0(x_1)=\sign(x_1)$ (so the $\eps\rightarrow0$ limit of $\varphi_\eps(x_1)$ is \eref{sign}). This has an invariant line $x_1=\eps c$, where $c$ is a constant such that $\varphi_\eps(\eps c)=0$, on which the dynamics is given by
\begin{equation}\label{slim0}
\sfrac{d\;}{dt}(x_1,x_2)=(0,1)\;.
\end{equation} 
This is illustrated in the far right of \fref{fig:eg1}. 
Now consider a piecewise smooth model of \eref{snon1}. In the limit $\eps=0$ equation \eref{snon1} becomes 
\begin{equation}\label{sdisc}
\sfrac{d\;}{dt}(x_1,x_2)=({\rm sign}(x_1),-1)\;. 
\end{equation}
Filippov's linear combination for \eref{sdisc}, using \eref{blowslidingf}, simplifies to 
\begin{equation}\label{eg1f}
\sfrac{d\;}{dt}(x_1,x_2)=(\lambda,-1)\;,
\end{equation}
where $\lambda$ is given by \eref{sign}. A sliding mode on the discontinuity surface $x_1=0$ must satisfy $\sfrac{d\;}{dt} x_1=0$ by \eref{blowslidingf}, so sliding modes exist for this system with $\lambda=0$, giving a sliding vector field $\sfrac{d\;}{dt}(x_1,x_2)=(0,-1)$, shown in the far left of \fref{fig:eg1}. Clearly this is not even close to being similar to the smooth system's dynamics in the far right portrait in \fref{fig:eg1}, given by \eref{slim0}. 
Moreover if we smooth this system, by replacing $\lambda$ with a sigmoid function (we can use $\varphi_\eps$ again for brevity), we obtain 
\begin{equation}\label{eg1sfs}
\sfrac{d\;}{dt}(x_1,x_2)=(\varphi_\eps(x_1),-1)\;,
\end{equation}
a system with an invariant line at some $x_1=c$, where $c$ is a constant such that $\varphi_\eps(x_1c)=0$, but on this line the vector field is again $\sfrac{d\;}{dt}(x_1,x_2)=(0,-1)$, inconsistent with \eref{slim0}; this is the second phase portrait in \fref{fig:eg1}. 

We obtained inconsistent sliding dynamics here, which persists under smoothing, because we did not keep the hidden term $\varphi_\eps^2$, neglecting is because is just unity for all $x_1\neq0$. If we take \eref{snon1} again, and carefully replace $\varphi_\eps(x_1)$ with $\lambda$ in the limit $\eps\rightarrow0$, we obtain 
\begin{equation}\label{snonfil}
\sfrac{d\;}{dt}(x_1,x_2)=(\lambda,1-2\lambda^2)\;.
\end{equation}
A sliding mode on the switching surface $x_1=0$ must satisfy $\sfrac{d\;}{dt} x_1=0$, and this has solutions for $\lambda=0$, which now gives a sliding vector field $\sfrac{d\;}{dt}(x_1,x_2)=(0,-1)$. This correctly captures the dynamics of the smooth system \eref{slim0} at the transition, and is illustrated by the third portrait in \fref{fig:eg1}. Moreover, obviously if we now smooth this by replacing $\lambda$ with a smooth sigmoid function, we regain \eref{snon1}, the far right portrait in \fref{fig:eg1}. 
\end{example}

In this example we obtained sliding dynamics that was not only quantitatively wrong compared to the system we were trying to model, but traveled in entirely the wrong direction. This happened because we neglected terms nonlinear in the switching term $\varphi_\eps$ or $\lambda$. 

The point is that if we start from the piecewise smooth model \eref{sdisc}, then either \eref{eg1f} or \eref{snonfil} are valid ways of continuing the system across $x_1=0$, amongst an infinity of other choices. In the following sections we will show how to express those possible systems, separate `linear' and `nonlinear' effects, and we shall highlight some more and less subtle effects they give rise to. 

A similar example to that above reveals that nonlinear terms can even affect whether or not the switching surface is crossed. 

\begin{example}
Consider a system 
\begin{equation}
\sfrac{d\;}{dt}(x_1,x_2)=(2\varphi_\eps^2(x_1)-1,1)\;.\label{snon2}
\end{equation} 
The idea is similar to the previous example so we just outline the calculations. This system has two invariant sets $x_1=\eps c_\pm$ such that $\varphi_\eps(\eps c_\pm)=\pm1/\sqrt2$, on which the vector field is $\sfrac{d\;}{dt}(x_1,x_2)=(0,1)$. In the limit $\eps=0$ we have
\begin{equation}
\sfrac{d\;}{dt}(x_1,x_2)=(1,1)\;,\label{eg2f}
\end{equation} 
from which the presence of a discontinuity is not even evident. Filippov theory therefore does not apply, and the system crosses through $x_1=0$. If we more carefully replace $\varphi_\eps$ by $\lambda$ in \eref{snon2}, we obtain 
\begin{equation}\label{eg2non}
\sfrac{d\;}{dt}(x_1,x_2)=(2\lambda^2-1,1)\;,
\end{equation} 
which has two sliding modes (where $\sfrac{d\;}{dt} x_1=0$) given by $\lambda=\pm1/\sqrt2$ implying $\sfrac{d\;}{dt}(x_1,x_2)=(0,1)$, equivalent to the invariant dynamics of the smooth system. Moreover, when this system is smoothed by replacing $\lambda$ with a sigmoid function like $\varphi_\eps$, we of course return to the original system \eref{snon2}. 
\end{example}
\begin{figure}[h!]\centering\includegraphics[width=0.45\textwidth]{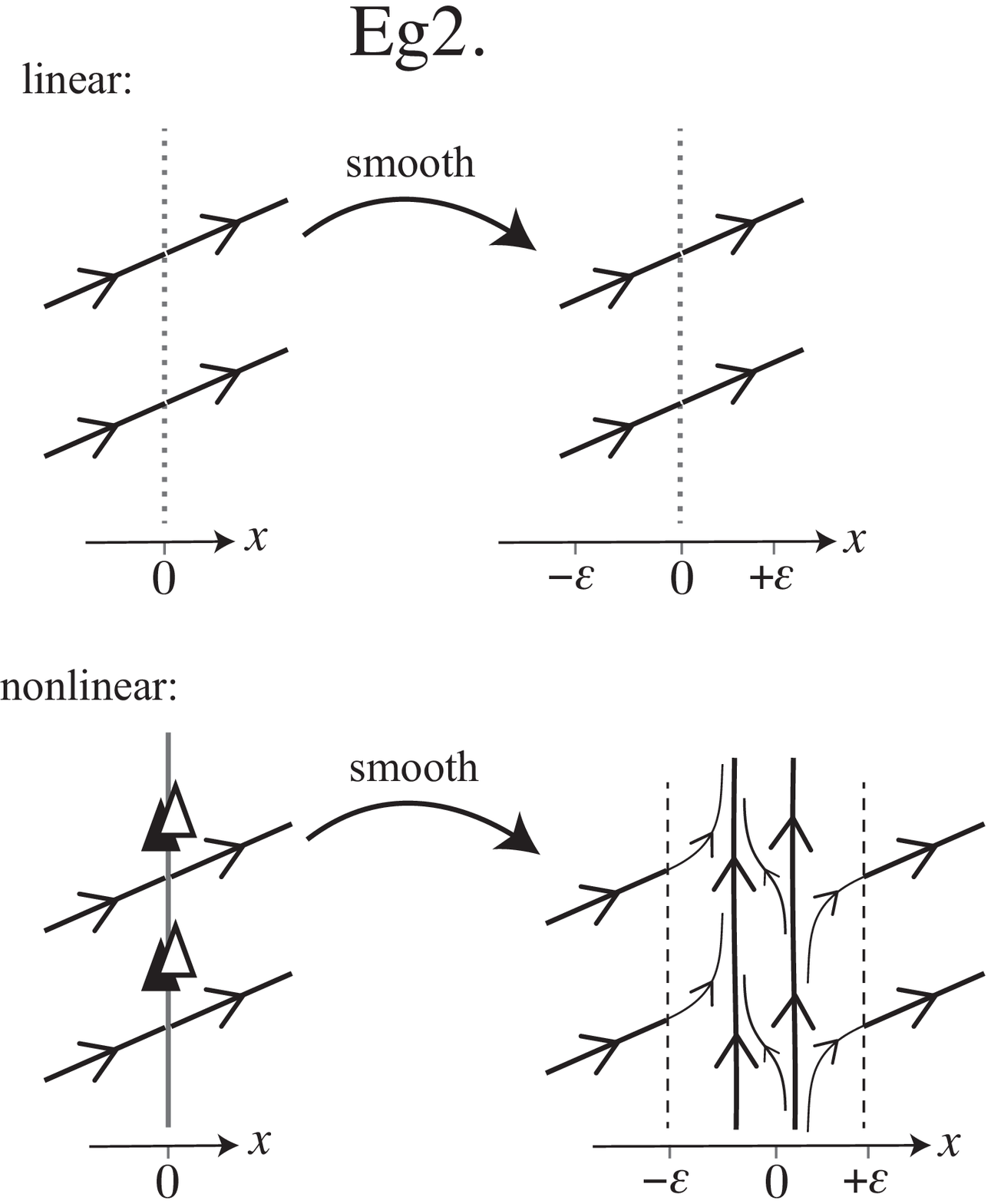}
\vspace{-0.3cm}\caption{\sf The system $\sfrac{d\;}{dt}(x_1,x_2)=(2({\rm sign}(x_1))^2-1,1)$ viewed ignoring the nonlinear dependence on $\sign(x_1)$ (left) or respecting it (right), in each case showing the discontinuous system and its smoothing. The phase portraits as shown correspond to equations \eref{eg2f}, \eref{eg2f}, \eref{eg2non}, \eref{snon2}, respectively. }\label{fig:eg2}\end{figure}

In this example the piecewise smooth model outside $x_1=0$ is given by \eref{snon2}, which is formed by patching together two systems for $x_1>0$ and $x_1<0$. Even though the vector fields appear continuous, their regions of definition are not. We may then assume the system {\it is} continuous, or consider \eref{snonfil} (and again an infinity of other choices) as a valid continuation of the system across $x_1=0$. In this case, the assumption of simple continuity turns out to be incorrect for modeling the system \eref{snon2} that we started from.

These are just trivial systems, and in the next section we will show how vast the field of possible regularizations is. Yet the ambiguities these examples highlight are extremely important to understanding the modeling of switches. Other examples can be much more subtle, as we shall see later. In \sref{sec:sigseries} we discuss how to handle nonlinearities to avoid ambiguity.

\section{The smoothing tautology -- keeping the hidden terms}\label{sec:smooth}


What is not clear in many works concerning smoothing is that there are many ways a system can be smoothed, and that they result in non-equivalent systems, any of which is equally valid as an approximation of the piecewise model \eref{ns1}. In this section we make this ambiguity explicit in the piecewise smooth system, as a `hidden' term in the vector field. 

Consider smooth systems expressible in the form
\begin{equation}\label{pert1}\begin{array}{lllll}
&\frac{dx}{dt}={f}\bb{{x};\varphi_\eps\!\bb{\sw({x})}}\;\qquad{\rm s.t.}\\\\
&\varphi_\eps\bb{\sw}\in\left\{\begin{array}{lll}\sign\bb{\sw}&\rm if&|\sw|>\eps\\\sq{-1,+1}&\rm if&|\sw|\le\eps\end{array}\right\}+\ord{\eps/\sw,\eps}
\end{array}\end{equation}
where the transition function $\varphi_\eps$ is at least continuous, and $\eps>0$ is small. The error term $\ord{\eps/\sw,\eps}$ ensures that we obtain \eref{ns1} as $\eps\rightarrow0$. We do not specify the form of $\varphi_\eps$ beyond this, for example it need not be monotonic, and need not be bounded to $\sq{-1,+1}$ at $\sw=0$. 


Let us moreover suppose $f$ to depend on $\varphi_\eps$ in an arbitrary manner, that is, it could be a polynomial, trigonometric, exponential or other function of the sigmoid. 
Nevertheless, if we take the limit $\eps\rightarrow0$ in \eref{pert1} we obtain \eref{ns1}. 
Does the dynamics of \eref{pert1} therefore tend to the dynamics of \eref{ns1} as described by \eref{ffil}?

Since we do not know precisely the form of $\varphi_\eps$ in the limit $\eps\rightarrow0$, let us introduce a simpler `reference' sigmoid for which we do. We define it as
\begin{equation}\label{leps}\begin{array}{lll}
&\Lambda_\eps(\sw)\in\left\{\begin{array}{lll}\sign(\sw)&\rm if&|\sw|>\eps\\\sq{-1,+1}&\rm if&|\sw|\le\eps\end{array}\right\}+\ord{\eps}\;,\\\\&{\rm s.t.}\qquad \Lambda_\eps'(\sw)>0\;\;\mbox{for}\;|\sw|<\eps\;,\end{array}
\end{equation}
and we shall use this to analyze the behaviour of $\varphi_\eps$ as $\eps$ tends to zero. 
Since $\Lambda_\eps(\sw)$ is differentiable and monotonic for $\eps>0$, at least for $|v|<\eps$, it has an inverse $V(\Lambda_\eps)$ such that $V(\Lambda_\eps(\sw))=\sw/\eps$. We can then define a new expression for $\varphi_\eps$ as a function of $\Lambda_\eps$, 
\begin{equation}
\Psi(\Lambda_\eps):=\varphi_\eps\bb{\eps V(\Lambda_\eps)}\;.
\end{equation}
Let us now use the reference sigmoid $\Lambda_\eps$ to form a linear combination of $f(x;\pm1)$,  
$$\op L(x,\Lambda_\eps(\sw))=\frac{1+\Lambda_\eps(\sw)}2f(x;+1)+\frac{1-\Lambda_\eps(\sw)}2f(x;-1)+. . .\;,$$
such that $\frac{dx}{dt}=\op L(x,\pm1)$ in the limit $|v|/\eps\rightarrow\infty$. This will not be exactly equal to \eref{pert1} for general $v$ and $\eps$, so taking the difference between this system and \eref{pert1}, we have
\begin{equation}\label{bigE}\begin{array}{lll}
&\frac{dx}{dt}=
\op L(x,\Lambda_\eps(\sw))+\op E\bb{x,\Lambda_\eps(\sw(x))}\;,\end{array}
\end{equation}
where $\op E\bb{x,\Lambda_\eps}=f\bb{x;\Psi(\Lambda_\eps)}-\op L(x,\pm1)$ 
vanishes when $\Lambda_\eps=\pm1$. 

Because the reference sigmoid $\Lambda_\eps$ varies monotonically with $\sw$ and is well-behaved as $\eps$ tends to zero, we can identify it with a new variable $\lambda$, then \eref{bigE} becomes
\begin{equation}\label{bigEl}
\frac{dx}{dt}=\frac{1+\lambda}2f(x;+1)+\frac{1-\lambda}2f(x;-1)+\op E\bb{x,\lambda}\;.
\end{equation}
The quantity $\lambda$ will vary dynamically according to $\sfrac{d\;}{dt}\lambda=\sfrac{d\lambda}{d\sw}\sfrac{d\;}{dt}\sw$, where $\sfrac{d\lambda}{d\sw}$ is strictly positive (which would not have been the case if we had identified $\lambda$ with $\varphi_\eps$). 

The expression \eref{bigEl} is now $\eps$ independent, and we may consider $\lambda$ to represent any $\varphi_\eps$, including the limiting function $\lambda=\phi_0(\sw)=\sign(v)$, with which \eref{bigEl} defines a piecewise smooth system consistent with \eref{ns1}, but notably more general than \eref{ffil}. 

Since $\op E$ vanishes for $\sw\neq0$, it is not fixed by \eref{ns1}. 
$\op E$ is the hidden part of the vector field, and by changing it we obtain arbitrarily many different non-equivalent dynamical systems, all consistent with \eref{ns1}.

There is of course no contradiction or paradox, nor indeed any surprise, in the fact that infinitely many different smooth systems expressible as \eref{bigE}, have the same limit \eref{ns1} as $\eps\rightarrow0$. This is merely due to the non-uniqueness of the singular limits by which smooth functions may tend towards discontinuities. What is important is that we can represent these different systems in the piecewise smooth limit by \eref{bigEl}. The system \eref{ffil}, therefore, is only the simplest member of a general class given by \eref{bigEl}.

%
%

The argument above began with \eref{pert1} as a prototype smooth system, but one may consider many other kinds. For instance consider a number of switches that activate at a common threshold $\sw=0$, with each switch having a different sigmoid-like behaviour $\varphi_1,\varphi_2,...$, so \eref{pert1} becomes
\begin{equation}\label{nreg}\begin{array}{ll}
&\frac{dx}{dt}={f}\bb{{x};\varphi_1\!\bb{\sw({x})},\varphi_2\!\bb{\sw({x})},...}\;\quad{\rm where}\\\\&\varphi_{\eps_ k}\bb{\sw}\in\left\{\begin{array}{lll}\sign\bb{\sw}&\rm if&|\sw|>\eps_k\\\sq{-1,+1}&\rm if&|\sw|\le\eps_k\end{array}\right\}+\ord{\eps_k/\sw}\;.\end{array}
\end{equation}
For example we might let
$$\begin{array}{rclcrcl}
\varphi_{\eps_1}(\sw)&=&\hf+\hf\frac{(\sw-1)^{1/\eps_1}}{1+(\sw-1)^{1/\eps_1}}\;,\\
\varphi_{\eps_2}(\sw)&=&\tanh(\sw/\eps_2)\;,\\
\varphi_{\eps_3}(\sw)&=&\erf(x/\eps_3)\;,\\
\varphi_{\eps_4}(\sw)&=&\frac 2\pi\arctan(\sw/\eps_4)\;,\quad etc.
\end{array}$$
We can again use a reference sigmoid $\Lambda_\eps$ with inverse $\sw=\eps V(\Lambda_\eps)$ to define a new expression for each $\varphi_{\eps_ k}$, 
\begin{equation}
\Psi_{k}(\Lambda_\eps):=\varphi_{\eps_ k}\bb{\eps V(\Lambda_\eps)}\;,
\end{equation}
and derive \eref{bigE} as above with $\op E\bb{x,\Lambda_\eps}=f\bb{x;\Psi_{1}(\Lambda_\eps),\Psi_{2}(\Lambda_\eps),...}-\frac{(1+\Lambda_\eps)f(x;+1)-(1-\Lambda_\eps)f(x;-1)}2$, which vanishes when $\Lambda_\eps=\pm1$, resulting again in the piecewise smooth system \eref{bigEl}. 

An interesting example arises when the sigmoids appear only linearly, but in a matrix form such as
\begin{equation}\label{mreg}
\frac{dx}{dt}=\frac{1+{\op M}(\sw(x))}2f(x;+1)+\frac{1+{\op M}(\sw(x))}2f(x;-1)\;,
\end{equation}
which we may call a {\it matrix regularization} of \eref{ns1}, where 
$${\op M}(\sw)=\bb{\begin{array}{cccccc}
\varphi_{\eps_1}\bb{\sw}&0&0&...\\
0&\varphi_{\eps_2}\bb{\sw}&0&...\\
0&0&\varphi_{\eps_ k}\bb{\sw}&...\\
\vdots&\vdots&\vdots&\ddots
\end{array}}$$
or even
$${\op M}(\sw)=\bb{\begin{array}{cccccc}
\varphi_{\eps_1}\bb{\sw}&\delta_{\eps_{12}}(\sw)&\delta_{\eps_{13}}(\sw)&...\\
\delta_{\eps_{21}}(\sw)&\varphi_{\eps_2}\bb{\sw}&\delta_{\eps_{23}}(\sw)&...\\
\delta_{\eps_{31}}(\sw)&\delta_{\eps_{32}}(\sw)&\varphi_{\eps_ k}\bb{\sw}&...\\
\vdots&\vdots&\vdots&\ddots
\end{array}}$$
where the off-diagonal elements are gaussian-like functions $\delta_{ij}(\sw)=\ord{\eps}$ for $|\sw|>\eps$ and $\delta_{ij}(\sw)\in[0,1]$ for $|\sw|\le\eps$. To handle these off-diagonal terms we need an extra step, because they are not sigmoids, nevertheless we can recast them again in terms of the reference sigmoid $\Lambda_\eps$ via a function $\Delta_{ij}:=\delta_{ij}\bb{\eps_{ij} V(\Lambda_{\eps_ k})}$, then to derive the corresponding hidden terms in $\op E$ is a simple exercise, resulting in \eref{bigEl} as before. 

The class of piecewise smooth dynamical systems that arise from such considerations is clearly vast, but all cases are encompassed in \eref{bigEl} with a suitable hidden term $\op E$~\footnote{\protect the obvious cases which are not in this class are those in which $\varphi_\eps$ depends on $x$ itself, not the threshold function $\sw$ alone, the crucial difference being that for these the transition cannot be reduced to a one dimensional system as in \eref{ldash}.}. 
%
Of course above we have written $\op E$ in terms of the systems \eref{pert1}, \eref{nreg}, or \eref{mreg}, expressions that are not useful in themselves (except that they show $\op E$ vanishes outside the switching surface). 
Some means to express the hidden term $\op E$ explicitly is now required and, where possible, to determine it from a system's behaviour. We broach this in the following section.

\section{Sigmoid asymptotics}\label{sec:sigseries}


Given a system $\sfrac{d\;}{dt}x=f$, and from only the knowledge that $f$ jumps between functional forms $f_+(x)$ and $f_-(x)$ as $\sw(x)$ changes sign, we wish to derive a complete model for the system. This will take the form $f=f(x;\lambda)$ in terms of a piecewise constant $\lambda$, or $f=f(x;\varphi_\eps(\sw))$ in terms of a sigmoid function $\varphi_\eps$. 

Below we will attempt formally to expand $f$ as a power series
\begin{equation}\label{asy0}
\frac{dx}{dt}=f(x;\lambda)=\sum_{n=0}^\infty \sa_n(x)\lambda^n(\sw(x))
\end{equation}
where for a piecewise smooth model we assume $\lambda=\sign(\sw)$, 
and for a smooth sigmoid model we replace $\lambda\mapsto\varphi_\eps(\sw)$, assuming $\varphi_\eps'(\sw)\ge0$ and
\begin{equation}\label{sigmoid}
\varphi_\eps(\sw)=\sign(\sw)+\ord\eps\qquad\mbox{for}\quad|\sw|>\eps\;.
\end{equation}

The possibility of approximating a function of several variables $(x,y)$ by a sum of functions of a single variable $\sw(x)$ has arisen particularly in the context of universal approximation by neural networks \cite{1989sigmoid,1998sigmoid,2004sigmoid}, where the sums of interest include sigmoid functions $\varphi(\sw_i(x))$ that may each have a different threshold $\sw_i=\hat k_i\cdot(x-x_i)$ and stiffness $|k_i|$ for some vectors $\cc{k_i,x_i}_{i=1,2,...}$ (falling therefore into our subclass of matrix regularizations in \eref{mreg}). Here we are interested in a single switching threshold, say $x_i=0$, and will be interested in the limit of infinite stiffness $|k_i|\to\infty$. An alternative is to use a series of powers of the sigmoid function $\varphi\bb{\sw}$, and the use of such series as solutions of nonlinear differential equations has been explored in \cite{tanh2004,tanh2004w} for $\varphi=\tanh$. 

The task now is to find the coefficients of the power series \eref{asy0}. 

\subsection{A sigmoid sum}

Comparing \eref{sigmoid} to \eref{ns1} in the states $\lambda=\pm1$, we have immediately that
\begin{equation}\label{pm0}
\begin{array}{lll}
&\displaystyle\sum_{n=0}^\infty \sa_n(x)=f_+(x)\qquad&\bb{\mbox{from }\sw>0}\;,\\
&\displaystyle\sum_{n=0}^\infty \sa_n(x)(-1)^n=f_- (x)\qquad&\bb{\mbox{from }\sw<0}\;.\end{array}
\end{equation}
The sum and difference of these give $ 2\sum_{n\;\rm even}\sa_n=f_++f_- $ and $ 2\sum_{n\;\rm odd}\sa_n=f_+-f_- $, allowing us to eliminate $\sa_0$ and $\sa_1$ to obtain 
\begin{eqnarray}\label{fEsum}
\frac{dx}{dt}=\hf(f_++f_- )+\hf(f_+-f_- )\lambda(\sw)+\op E(x;\lambda)\;,
\end{eqnarray}
where
\begin{eqnarray}\label{Eandg}
\op E(x;\lambda)&=&\sum_{n=1}^\infty \bb{\sa_{2n}+\lambda(\sw) \sa_{2n+1}}\bb{\varphi^{2n}(\sw)-1}\nonumber\\
&=&(\lambda^2-1)\sum_{n=1}^\infty\sum_{j=0}^{n-1}\sq{\sa_{2n}(x)+ \lambda\sa_{2n+1}(x)}\lambda^{2j}\nonumber\\
&:=&(\lambda^2-1)\;g(x,\lambda)\;,
\end{eqnarray}
simply factoring out $\lambda^2-1$ in the second line, and defining the double sum as a function $g$. Note that while $\op E$ vanishes for $\lambda=\pm1$, i.e. outside the switching surface, the function $g$ can take any value. Hence we obtain exactly \eref{bigEl}, and a formula for $\op E$. 

In the piecewise smooth case, the remaining coefficients $\alpha_{n\ge2}$ that make up the hidden term $g$ may be determined if we have any information about $f$ at $\sw=0$. For example, say it is known that for $\lambda=0$ we have $\left.f(x;0)\right|_{\sw=0}=r(x)$, then $\alpha_{2}=\hf(f_++f_-)-r(x)$, and successive coefficients can be eliminated with knowledge of derivatives of $f$ with respect to $\lambda$. If these quantities cannot be derermined by direct observations of a system, we may propose forms of $g$ to fit dynamical observations or physical laws, depending on the modeling problem. In the smooth case, where $\lambda$ is replaced by a differentiable sigmoid function $\phi_\eps$, then more can be found by considering the asymptotics of the transition.

\subsection{Asymptotics of the series}


Here we consider how \eref{fEsum}-\eref{Eandg} may be applied to a system where the transition between $f_\pm$ is smooth. 

In a system that makes a sudden transition between behaviours $\frac{dx}{dt}=f_+(x)$ for large positive $\eta=\sw/\eps$, and $\frac{dx}{dt}=f_- (x)$ for large negative $\eta=\sw/\eps$, let us assume that $f_+$ and $f_- $ are actually steady states $f(x;y_+)$ and $f(x;y_- )$ of a larger bistable system in $x$ and $y$, 
such that
\begin{equation}\begin{array}{rcl}
y&\rightarrow&\left\{\begin{array}{lll}y_+&\rm as&\eta\rightarrow+\infty\;,\\y_- &\rm as&\eta\rightarrow-\infty\;,\end{array}\right.\\\\\Rightarrow
\quad\frac{dx}{dt}&\rightarrow&\left\{\begin{array}{lll}f_+&\rm as&\eta\rightarrow+\infty\;,\\f_- &\rm as&\eta\rightarrow-\infty\;.\end{array}\right.
\end{array}
\end{equation}
The vector $y$ may represent microscopic properties or external forces independent of $x$, which vary significantly only where their value jumps. 
The scalar $\eta$ may merely be a fast timescale $t/\eps$ for small $\eps$, or else some function $\eta(x,\eps)$ of the internal state
. Without specific knowledge of the $y$ system, we seek to model the effect on the system $\frac{dx}{dt}=f(x)$. 

If we could linearize with respect to $y$ near the steady state $f_+=f(x;y_+)$ for large negative $\eta$, 
treating $\eta$ as a time-like variable, then $f$ should diverge from $f_+$ along a direction $g_+$ in phase space with a time constant $\kappa$ (the Lyapunov exponent). The resulting approximation $df/d\eta\sim\kappa\bb{f-f_+}$ has a solution of the form $f\sim f_++g_+ e^{\kappa\eta}\sum_{n=0}^\infty \sb_n\eta^{-n}$ for some coefficients $\sb_n$. 
This is easily extended if we cannot linearize. For example, if $f=f_+$ is a degenerate steady state of the $y$ system, the same asymptotic formula applies with $\kappa=0$ if $f=f_+$ is a saddle-node, where $df/d\eta\sim|f-f_+|(f-f_+)/\sb_1$ for some constant $\sb_1$. The same formula again holds with the exponent replaced by $-\kappa|\eta|^p$ near a non-hyperbolic steady state of the form $df/d\eta\sim\kappa|\eta|^{p-1}(f-f_+)$. 

Based on the typical divergence from a linear or nonlinear steady state, therefore, near either of the states $\frac{dx}{dt}=f_+$ or $\frac{dx}{dt}=f_- $, the system behaves asymptotically like 
\begin{equation}\label{asy1}
\frac{dx}{dt}\sim f_i(x)+g_i(x)e^{-\kappa_i|\sw/\eps|^{p_i}}\sum_{n=0}^\infty \sb_n^{(i)}\!(x)(\eps/\sw)^{n}
\end{equation}
as $\sw/\eps\rightarrow-\infty$ for $i=+$ or $\sw/\eps\rightarrow+\infty$ for $i=-$, 
for positive constants $\kappa_{i},p_{i},\sb_n^{(i)}$. The vector-valued functions $g_i(x)$ describe the trajectories in function space by which $\frac{dx}{dt}$ departs from each of the $f_i$. 
Provided $\sb_0^{(i)}=0$ if $\kappa_i=0$, this ensures that the second term vanishes for $\sw/\eps\rightarrow\pm\infty$. 

If we had the exact differential equations for $f(x;y)$, then we would now apply asymptotic matching (see e.g. \cite{bo99}) to unify the two approximations in \eref{asy1}. When no such differential equations are known, we may seek to represent $f$ by a series approximation {\it in $x$ only}, since $y$ is assumed to be beyond ready observability. 
The terms of the series must be functions that give the correct asymptotic behaviour for positive and negative $\sw$, and this can be achieved with sigmoid functions $\varphi$ that vary monotonically from $-1$ to $+1$, switching value at $\eta(x)=0$. 
Thus we take the series \eref{asy0} with \eref{sigmoid}. 

Finding the first two coefficients is similar to the piecewise smooth case. 
Comparing \eref{asy0} to \eref{asy1} at the steady states $f_+$ and $f_- $ we have that $\alpha_0=\hf(f_++f_-)+ \sum\sa_{2n}$ and $\alpha_1=\hf(f_+-f_-)+ \sum\sa_{2n+1}$, so
\begin{eqnarray}
\frac{dx}{dt}
=\frac{f_++f_- }2+\frac{f_+-f_- }2\varphi_\eps(\sw)+\op E(x,\varphi_\eps(\sw))
\end{eqnarray}
where $\op E(x,\varphi_\eps)=(\varphi_\eps^2-1)g(x,\varphi_\eps)$, with $\op E$ and $g$ being the same functions as in the piecewise smooth case \eref{Eandg}. 

The remaining coefficients gathered together in $g$ can be found by assuming that $\varphi_\eps$ has an asymptotic expansion of the form 
\begin{equation}\label{sigasy}
\varphi_\eps(\sw)\sim\sign(\sw)\bb{1+e^{-\kappa_\sc|\sw/\eps|^{p_\sc}}\sum_{n=0}^\infty \sc_n(\eps/\sw)^{n}}\;,
\end{equation}
consistent with \eref{sigmoid}. For example, the Hill functions \cite{hill} behave as $Z(x)=\frac{x^{1/\eps}}{x^{1/\eps}+\theta^{1/\eps}}\sim\hf+\sign(x-\theta)\bb{\hf-e^{-\abs{\log(x/\theta)}/\eps}+...}$  for $x>0$, and we let $\sw=x-\theta$, $\varphi_\eps(\sw)=2Z(\theta e^{\sw/\eps})-1$, and other common sigmoids include $\varphi_\eps(\sw)=\tanh(\sw/\eps)\sim\sign(\sw)\bb{1-2\eps^{-2|\sw|/\eps}+...}$, $\varphi_\eps(\sw)=\erf(v/\eps)\sim\sign(\sw)-e^{-\sw^2/\eps^2}((\eps/\sw)-\hf(\eps/\sw)^{3}+...)/\sqrt\pi$.

To give an example of the coefficients $\sa_{n\ge2}$, if we assume some symmetry in the asymptotic limit such that $\kappa_+=\kappa_- $ and $p_+=p_- $, then it makes sense to choose $\kappa_\sc=\kappa_+=\kappa_- :=\kappa$ and $p_\sc=p_+=p_- :=p$, then comparing terms of order $e^{-k|\eta|^{p}}$ allows us to find $\sa_2$ and $\sa_3$ as
\begin{equation}
\begin{array}{lll}
\sa_{2}&=&\frac1{4\sc_0}\bb{g_+\sb_0^{(a)}+g_- \sb_0^{(b)}}-\sum_{m=2}^\infty m\sa_{2m}\;,\\
\sa_{3}&=&\frac1{4\sc_0}\bb{g_+\sb_0^{(a)}-g_- \sb_0^{(b)}}\\&&\qquad\quad-\sfrac14(f_+-f_- )-\sum_{m=2}^\infty m\sa_{2m+1}\;,\end{array}
\end{equation}
and so on, recursively for each coefficient $\sa_n$.

\section{Example: switching asymptotics in a control circuit}\label{sec:circuit}

The most obvious kind of sigmoidal behaviour in physics is a mechanical or electronic binary relay. Consider a typical electrical model based on a DC-DC converter, with a battery voltage $V_b $ and inductor current $I$. A connecting circuit, with a capacitance $C$ and resistance $R$ across which the voltage is $V$, is disconnected if $V$ exceeds a value $V_b $. The empirical model of the on and off positions is
\begin{equation}\label{onoff}\begin{array}{rcl}
\mbox{``on'':}&&\;\;\left\{\begin{array}{rcl}
L\frac{d\;}{dt}I&=&V_0 -V\;,\\
RC\frac{d\;}{dt}V&=&IR-V\;,
\end{array}\right.\\\\
\mbox{``off'':}&&\;\;\left\{\begin{array}{rcl}
L\frac{d\;}{dt}I&=&V_0 \;,\\
RC\frac{d\;}{dt}V&=&-V\;.
\end{array}\right.\end{array}
\end{equation}
%
\begin{figure}[h!]\centering\includegraphics[width=0.4\textwidth]{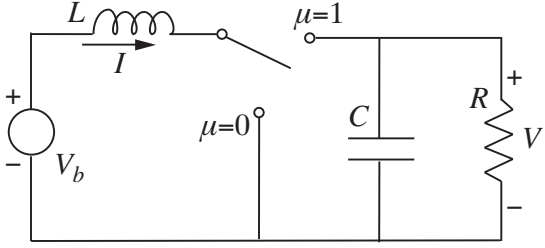}
\vspace{-0.3cm}\caption{\sf A basic switching circuit. }\label{fig:}\end{figure}
%

What then happens during the transition between these two? 
The simplest model we might write down, and the canoncial Utkin/Filippov approach \cite{f88,u77}, expresses \eref{onoff} as
\begin{equation}\label{circuitstep}\begin{array}{ll}
&\bb{L\sfrac{d\;}{dt}I,RC\sfrac{d\;}{dt}V}=f=\bb{V_0,-V}+\bb{-V,IR}\mu\\\\
&\qquad\qquad\qquad\quad{\rm where}\qquad\mu=\step(V_b-V)\;.
\end{array}\end{equation}
The switch $\mu$ is related to the multiplier $\lambda$ in the previous section by $\mu=\hf+\hf\lambda$, we use $\mu$ here as it is more common in this application.

Let us then consider a more general class of systems
\begin{equation}\label{circuitmodel}\begin{array}{ll}
&\bb{L\sfrac{d\;}{dt}I,RC\sfrac{d\;}{dt}V}=f=\bb{V_0,-V}+\bb{-V\mu,IRp(\mu)}\\\\
&\qquad\qquad\qquad\quad{\rm where}\qquad\mu=\step(V_b-V)\;,\end{array}
\end{equation}
such that $p(\mu)=\step(V_b-V)$. 
(This class is introduced for demonstration, and it may be possible to motivate more experimentally relevant classes from the precise physics of a given device). For the sake of an explicit example we shall consider $p(\mu)=\mu-\sigma(1-\mu)\mu$ below, 
so that the original system \eref{circuitstep} corresponds to $\sigma=0$. 
We shall ask to what extent we can distinguish between systems with different $\sigma$. 


The smooth function replacing the step should be monotonic and continuous, say $\mu=\varphi_\eps(V_b-V)$ defined in terms of a stiffness parameter $\eps$, such that $\varphi_0(V_b-V)=\step(V_b-V)$. The resulting system
\begin{equation}\label{circuitmu}\begin{array}{ll}
&\bb{L\sfrac{d\;}{dt}I,RC\sfrac{d\;}{dt}V}=f=\bb{V_0,-V}+\bb{-V\mu,IRp(\mu)}\\\\
&\qquad\qquad\qquad\quad{\rm where}\qquad\mu=\varphi_\eps(V_b-V)\;,\end{array}
\end{equation}
becomes \eref{circuitmodel} as $\eps\rightarrow0$, i.e. this is an $\eps$-perturbation of \eref{circuitmodel}.  

The system \eref{circuitmu} has a small deviation from the steady states $\mu=0$ and $\mu=1$ for $|\sw|>\eps$, described by the derivative with respect to $V/\eps$, 
\begin{equation}\label{circd}
\frac{\partial \;}{\partial V/\eps}f=\bb{-V,IRp'(\mu)} \frac{d\mu}{dV/\eps}+\ord\eps\;.
\end{equation}
The $p'$ term means this deviation will be $\sigma$ dependent. The term $\frac{d\mu}{dV/\eps}$ means the deviation will be of order $\ord{|V_b-V|/\eps}$ or smaller, as discussed in the previous section leading to \eref{asy1}. 
%
Take for example $\varphi_\eps$ to be an $\arctan$ function, specifically 
\begin{eqnarray*}
\varphi_\eps(V_b-V)&=&\hf+\sfrac1\pi\arctan\bb{\sfrac{V_b-V}\eps}\\&=&\step(V_b-V)-\sfrac{\eps}{\pi(V_b-V)}+\ord{\sfrac{\eps^2}{(V_b-V)^3}}
\end{eqnarray*}
then
\begin{eqnarray*}
\frac{d\mu}{dV/\eps}&=&\frac{d\varphi_\eps(V_b-V)}{dV/\eps}=\frac{\pi^{-1}}{1+|V_b-V|^2/\eps^2}\\&=&\frac{\eps^2}{2|V_b-V|^2}+\ord{\sfrac{\eps^3}{(V_b-V)^3}}\;,
\end{eqnarray*}
hence 
\begin{equation}\label{circd}
\frac{\partial \;}{\partial V/\eps}f=\bb{-V,IRp'(\mu)}\sfrac{\eps^2}{2|V_b-V|^2}+\ord{\sfrac{\eps^3}{(V_b-V)^3},\eps}\;.
\end{equation}

The difference between models with different $\sigma$ is therefore very small, skrinking with $\eps$ and with $\frac{\eps^2}{|V_b-V|^2}$ for $V$ away from $V_b$ (and alternative sigmoid functions may give deviations exponentially rather than polynomially small in $|V_b-V|/\eps$). Observing the deviation directly, therefore, would be a challenge in practice. Nevertheless we will show that $\sigma$ has significant and observable effects on the stability of the circuit no matter how small $\eps$ is (i.e. how close the system is to discontinuous), because the term $\frac{\eps^2}{|V_b-V|^2}$ permits significant deviations near $V=V_b$. 
So we now turn to the circuit's behaviour. 



\begin{figure}[h!]\centering\includegraphics[width=0.4\textwidth]{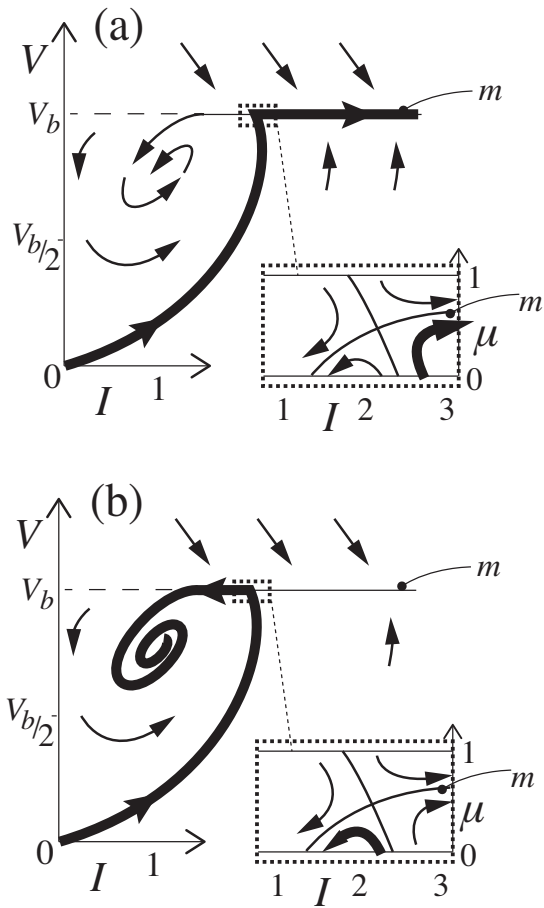}
\vspace{-0.3cm}\caption{\sf Phase portrait of the current, showing the piecewise smooth system (main picture), and a blow up of the switching surface into the layer $\mu\in[0,1]$ revealing a saddle. The trajectory from $I=V=0$ is shown, simulating values $L=5$, $V_0=5$, $V_b=6$, $RC=5/2$, $R=15/4$, and (a) $\sigma=0$, (b) $\sigma=1/2$. }\label{fig:smode}\end{figure}

With the switch in the ``on'' position the state spirals in towards an attracting focus at $I=V_0/R$, $V=V_0 $, as shown in \fref{fig:smode}. In the ``off'' position the system is driven towards the switching surface $V=V_b$. 
We then need a law describing how the systems evolves at the switch. 

We can use the smooth system to derive how $\mu$ jumps between $0$ and $1$ across the switching surface. Differentiating $\mu=\varphi_\eps(V_b-V)$ with respect to $t$ gives
\begin{equation}
\eps\frac{d}{dt}\mu=\eps\frac{d V}{dt}\frac{d}{dV}\mu\quad\Rightarrow\quad \tilde\eps(\eps,V)\frac{d}{dt}\mu=\sfrac{d\;}{dt} V\;,
\end{equation}
where $\tilde\eps=\eps/\frac{d\mu}{dV/\eps}$ satisfies $\eps<\tilde\eps<4\eps$ inside the transition region $|V-V_b|<\eps$, hence $\tilde\eps=\ord\eps$ is small around the switch. We can therefore introduce a fast timescale $\tau=t/\tilde\eps$ in terms of which $\frac{d\;}{d\tau}=\frac{dV}{dt}=f\cdot\nabla V$. The dependence on $\eps$ is supressed (but for $\eps=0$ the $\tau$ timescale is instantaneous compared to $t$), and we can take the limit $\eps\rightarrow0$ to obtain the piecewise smooth system, in which the transition region becomes the switching surface $V=V_b$, on which we have
\begin{eqnarray}\label{circuitblow}
\bb{L\frac{d\;}{dt}I,RC\frac{d\;}{d\tau}\mu}&=&\;\bb{V_0,-V_b}+\bb{-V_b\mu,IRp(\mu)}\;.\;\;
\end{eqnarray}

Thus \eref{circuitblow} prescribes the dynamics on $(I,\mu)$ at a point on the switching surface $V=V_b$. Let us take the function $p(\mu)$ as 
\begin{equation}
p(\mu)=\mu-\sigma(1-\mu)\mu
\end{equation} 
as an example, and assume $|\sigma|<1$. Then for $IR<V_b$ we have $d\mu/d\tau<0$, implying that the flow crosses directly from $V>V_b$ to $V<V_b$. For $IR>V_b$ instead, the interval or {\it layer} of the transition, $\mu\in[0,1]$, contains an attractor in the form of a saddlepoint at $\mu=\frac{V_0}{V_b}$, $IR=\frac{V_b^3}{V_0V_b-\sigma V_0(V_b+V_0)}$. The dynamics outside the layer are shown in the main part of \fref{fig:smode}, with a blow up (inset) of the dynamics on $\mu\in[0,1]$ inside the switching surface $V=V_b$, shown near the saddle. 

Solutions that enter the layer $\mu\in[0,1]$ will be attracted to one of the two unstable manifolds of the saddle, either traveling toward increasing $I$ with $V$ remaining on the switching threshold (shown in \fref{fig:smode}(a)), or traveling toward decreasing $I$, before detaching from the switching surface at $IR=V_b$ and spiralling in towards the focus in the ``on'' position (shown in \fref{fig:smode}(b)). 

The critical effect of $\sigma$ in this system is to shift the current $I=\frac{V_b^3/R}{V_0V_b-\sigma V_0(V_b+V_0)}$ at which the saddlepoint occurs on the switching surface. \Fref{fig:smode} shows the behaviour of the system from an initial condition $(I,V)=(0,0)$ for two different values of $\sigma$. 


If we set $\sigma=0$, then from $I=V=0$ the system trajectory enters the sliding mode to the right of the saddle, and the current grows unboundedly, while the voltage remains fixed at $V=V_b$. In this case the switch depends linearly on $\mu$, and the sliding mode is simply as would be obtained from Filippov/Utkin's standard method \cite{f88,u77}, where the saddlepoint at 
$\mu=\frac{V_0}{V_b}$, $IR=\frac{V_b^2}{V_0}$ is usually called a {\it psuedoequilibrium}. 

For $\sigma$ nonzero, instead, the saddlepoint moves to a higher current, and from the initial condition $I=V=0$ the system hits the sliding mode to the left of the saddle, the current falls until the sliding mode terminates, and the system then spirals in towards the focus. The difference in the dynamics at $V=V_b$ has to do with the fact that, with $\sigma\neq0$, the system is nonlinear in the switching parameter $\mu$. 


The effect of nonlinearity in the voltage response to the switch here is relatively minor -- shifting the saddlepoint -- but with stark effects on the system stability. The effects can easily be more interesting, as we show with an abstract but novel example in \sref{sec:ueda}. First we shall generalize the methods above for analyzing the piecewise smooth system, without having to concern ourself with the asymptotics of smoothing, yet without losing the distinction between systems that behave differently near the switch. 


\section{Nonsmooth dynamics -- the nonlinear approach}\label{sec:non}

We can generalize the methods used in the previous section with a slight development of methods proposed in \cite{f88,j13error}.

We have shown in \sref{sec:sigseries} that we can express the system \eref{ns1} as
\begin{equation}\label{fgen}
\frac{dx}{dt}=\frac{f_+(x)+ f_- (x)}2+\frac{f_+(x)- f_-(x) }2\lambda+\op E(x;\lambda)
\end{equation}
where  
\begin{equation}\label{Evanish}
\lambda=\sign\sw\qquad{\rm and}\qquad \sw\op E(x;\lambda)=0\;.
\end{equation}
The term $\op E$ is the `hidden' part of the vector field, the condition $\sw\op E=0$ meaning that this vanishes outside the switching surface. 

At $\sw=0$ the system \eref{ns1} with \eref{sign} gives a differential inclusion, but (as we found in the last section) it turns out that the value of $\lambda$ can be fixed in most cases by dynamical arguments. For $\op E\equiv0$ this is done using Filippov's sliding modes \cite{f88}, and in \cite{j13error} these ideas were extended to the case when $\op E$ is nonzero, i.e. when \eref{fgen} depends nonlinearly on $\lambda$. We shall review the method briefly. 

We introduce a transition timescale $\tau$ which is instantaneous in the original timescale $t$ (so we may think of this as $\tau=\lim_{\eps\to0}t/\eps$ for some $\eps>0$, as we obtained in the previous section). Since $\lambda$ switches between $\pm1$ as $\sw$ changes sign, its variation is driven by the component of the vector field along $\nabla\sw$ normal to the switching surface, so we define the dynamics of $\lambda$ as
\begin{equation}\label{ldash}
\frac{d\;}{d\tau}\lambda=f(x;\lambda)\cdot\nabla\sw\;\qquad{\rm on}\quad\sw=0\;.
\end{equation}
This is the generalization of the expression we obtained in the circuit model in \sref{sec:circuit}, and can be more generally derived from a smoothing of \eref{ns1} using arguments of asymptotics and geometric singular perturbation theory; see  \cite{j13error,j15hidden,j15sens2fold}. 

The expression \eref{fgen} with \eref{ldash} is now sufficient to prescibe the dynamics of the piecewise smooth system fully. Combining them at the switching surface, taking coordinates $x=(x_1,x_2,...,x_n)$ in which $x_1=\sw(x)$, and letting $f=(f_1,f_2,...,f_n)$, we have the two timescale system
\begin{equation}\label{blow}
\left.\begin{array}{rcl}\frac{d\;}{d\tau}\lambda&=&f_1({x};\lambda)\\
\frac{d\;}{dt}\bb{x_2,...,x_n}&=&\bb{f_2({x};\lambdab),...,f_n({x};\lambdab)}\end{array}\right\}
\end{equation}
on $x_1=0$. 
When a solution reaches the discontinuity, the system \eref{blow} facilitates either the transition from one side of the switch to the other, or, if equilibria of the $\tau$ timescale subsystem exist (i.e. where $\frac{d\;}{d\tau}\lambda=0$), then $\lambda$ collapses to a value $\lambda^s$ given by
\begin{equation}\label{blowsliding}
\left.\begin{array}{rcl}0&=&f_1({x};\lambda^s)\\
\frac{d\;}{dt}\bb{x_2,...,x_n}&=&\bb{f_2({x};\lambdab^s),...,f_n({x};\lambdab^s)}\end{array}\right\}
\end{equation}
on $x_1=0$. 
These equilibria, if they exist, represent {\it sliding modes}, whereby the dynamics sticks to the switching threshold $\sw=0$ and evolves ({\it slides}) along it. We see that these are precisely the sliding modes defined in \eref{blowslidingf} by Filippov, except that now with \eref{fgen} and \eref{blow} they apply also to systems with nonlinear dependence on $\lambda$.

These are the basic elements for solving a system with nonlinear switching. 
We have seen in three examples how nonlinear terms can alter the dynamics. Nonlinearity can lread to much more interesting dynamical phenomena, an example of which is given below. 

\section{Example: hidden Duffing oscillator}\label{sec:ueda}

\begin{figure}[h!]\centering\includegraphics[width=0.35\textwidth]{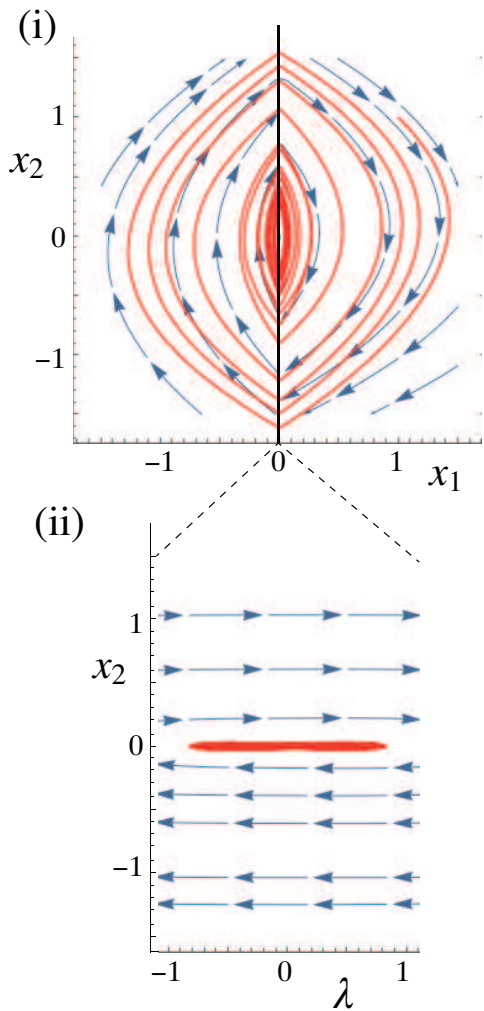}
\vspace{-0.3cm}\caption{\small\sf The flow of: (i) the fused focus \eref{ueda}, and (ii) its switching layer \eref{uedablow}, shown at a fixed time, for constants $c=0.1$, $b=0.05$, $a=0.15$. An orbit is simulated over a time $t$ from 0 to 500. }\label{fig:ueda}\end{figure}

Consider the system
\begin{equation}\label{ueda}
\sfrac{d\;}{dt}( x_1, x_2)=\bb{x_2-cx_1,\;-\lambda^3-bx_2+a\cos t}\;,
\end{equation}
where $\lambda=\sign(x_1)$, 
for constants $a,c,b$. The dynamics is that of a {\it fused focus}, that is, a folding of the dynamics either side of the switch that creates focus-like attraction towards $x_1=x_2=0$; this is illustrated in \fref{fig:ueda}(i).

For $x_1\neq0$ this is indistinguishable from the linear switching system
\begin{equation}\label{uedalin}
\sfrac{d\;}{dt}( x_1, x_2)=\bb{x_2-cx_1,\;-\lambda-bx_2+a\cos t}\;,
\end{equation}
where $\lambda=\sign(x_1)$, 
which is well described by Filippov/Utkin sliding mode theory \cite{f88,u77,krg03}. 

Both \eref{ueda} and \eref{uedalin} are clearly models consistent with the system
\begin{equation}\label{uedans}
\sfrac{d\;}{dt}( x_1, x_2)=\left\{\begin{array}{lll}\bb{x_2-cx_1,\;-1-bx_2+a\cos t}&\rm if&x_1>0\;,\\\bb{x_2-cx_1,\;+1-bx_2+a\cos t}&\rm if&x_1<0\;.\end{array}\right.
\end{equation}

Let us apply the methods from the previous section to \eref{ueda} to see the effect of the $\lambda^3$ term in \eref{ueda}. Substituting \eref{ueda} into \eref{blow} to find the dynamics on $x_1=0$ gives
\begin{equation}\label{uedablow}
(\sfrac{d\;}{d\tau}\lambda,\sfrac{d\;}{dt} x_2)=\bb{x_2,\;-\lambda^3-bx_2+a\cos t}\qquad{\rm on}\quad x_1=0\;,
\end{equation}
where $\tau$ is an instantaneous transition timescale. Only the attracting point $x_1=x_2=0$ satisfies the equation \eref{blowsliding}, so sliding modes are confined to the origin and would seem to be of little interest. However, in the variables $(\lambda,x_2)$ the system \eref{uedablow} is a form of Duffing oscillator (see e.g. \cite{gh02}), and this can lead to rich dynamics inside the transition layer $\lambda\in[-1,+1]$. 

Simulations indeed reveal nontrivial dynamics on $\lambda$ and $x_2$ for $\lambda\in[-1,+1]$. One simulation is shown in \fref{fig:ueda}(ii), showing the switching surface $x_1=0$ `blown up' into a switching layer $\lambda\in[-1,+1]$ on $x_1=0$. The variable $x_2$ varies only slightly, while $\lambda$ oscillates in an irregular pattern between values of around $|\lambda|\le a^{1/3}$. 

The structure of the orbit in \fref{fig:ueda}(ii) is not visible, and may appear to be of little significance anyway, since it is negligible in the state space of $(x_1,x_2)$ in \fref{fig:ueda}(i). The switching paramater $\lambda$ may, however, have a noticable effect on the system, for example if it affects another variable. Let us introduce a third variable satisfying some simple equation involving the switch, say
\begin{equation}\label{uedaz}
\mu\sfrac{d\;}{dt} x_3=\lambda-x_3\;,
\end{equation}
where $\mu$ is some constant, then $x_3$ will be simply $\pm1$ outside the switching surface, but inside the surface it will track the dynamics of $\lambda$. We choose $\mu$ small enough in our simulations so that $x_3$ and $\lambda$ are indistinguishable.

\begin{figure}[h!]\centering\includegraphics[width=0.45\textwidth]{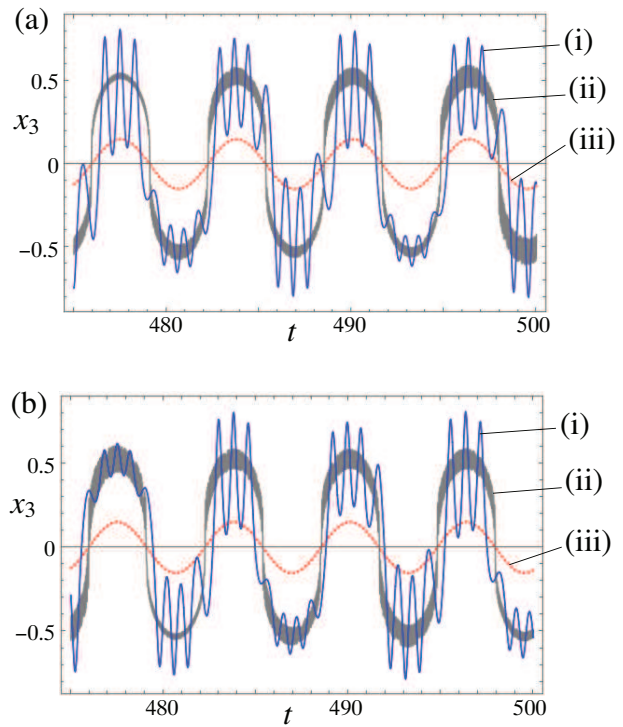}
\vspace{-0.3cm}\caption{\small\sf Simulations of hidden dynamics, plotting $z(t)$ (or $\varphi(x(t)/\eps)$ for: (i) the nonlinear system with $\eps=10^{-2}$, (ii) the nonlinear system with $\eps=10^{5}$, and (i) the linear system the $\eps=10^{-1}$. Panel (a) is a simulation of the switching layer system \eref{uedablow}, and panel (b) is a simulation of the full system \eref{ueda} smoothed out using $\lambda=\varphi_\eps(x_1)$. All simulations are with $c=0.1$, $b=0.05$, $a=0.15$.}\label{fig:uedag}\end{figure}

After discarding transients, $x_3(t)$ is seen to have oscillatory, complex, and chaotic dynamics for different parameters, including Ueda's chaotic attractor \cite{ueda92}, as expected from a Duffing oscillator. 
\Fref{fig:uedag}(a) shows a plot of $x_3(t)$. The speed of the timescale $\tau$ matters here, so we let $\tau=t/\eps$, the ideal instantaneous switch corresponding to the limit $\eps\rightarrow0$. The curves (i) and (ii) show $x_3(t)$ (or equivalently the value of $\lambda$) for $\eps=10^{-2}$ and $\eps=10^{-5}$, respectively. The oscillations become more regular as $\eps$ tends to zero, with the amplitude of this regular oscillation tending to around $a^{1/3}$. 

Contrast this with the curve (iii), which shows the graph of $x_3(t)$ for the linear model \eref{uedalin}, obtained by replacing $\lambda^3$ with $\lambda$ in \eref{uedablow}, and gives instead always a steady oscillation (which does not vary appreciably with $\eps$), but with an amplude $a$, rather than the $a^{1/3}$ of the nonlinear system. 

Throughout this paper we have tried to show that smoothing is a tautologous process, that is, it cannot be used to justify one way or another of resolving a discontinuous system, say taking \eref{ueda} or \eref{uedalin} to resolve \eref{uedans}, and that systems must be distinguished in the expression of the discontinuous systems themselves. We demonstrate this finally by smoothing out the systems \eref{ueda} and \eref{uedalin} directly, replacing $\lambda$ in those expressions with
\begin{equation}\label{uphi}
\varphi_\eps(x_1)=\left\{\begin{array}{lll}\sign(x_1)&\rm if&|x_1|\ge\eps\;,\\x_1/\eps&\rm if&|x_1|<\eps\;,\\\end{array}\right.
\end{equation}
and we simulate in \fref{fig:uedag}(b) the orbits corresponding to those in \fref{fig:uedag}(a). The qualitative dynamics is similar, with only minor quantitative differences in the precise form of smaller oscillations, confirming that smoothing only preserves the differences between the linear and nonlinear models. Simulations using alternative smoothing functions such as $\varphi_\eps(x_1)=\sfrac2\pi\arctan(x_1/\eps)$ or $\varphi_\eps(x_1)=\tanh(x_1/\eps)$ yield results with no significant difference to those using \eref{uphi}.

\section{Closing Remarks}\label{sec:conc}

The traditional use of Filippov's theories to tackle piecewise smooth systems neglects the possibility of nonlinear dependence on $\lambda$; (we should remark that Filippov's work \cite{f88} sets out results for much more general differential inclusions, but most of his detailed theory and most modern theory derived from it assume linear dependence). In many cases the assumption of linear dependence on $\lambda$ will be sufficient, but nonlinear dependence may have to be introduced if a linear model's dynamics that is inconsistent with experiment, if we smooth the system out, or, something we have not discussed here, if singularities arise that exhibit degeneracy inside the switching surface \cite{j15sens2fold}. 

Great intricacy is to be found in not only the global, but the {\it local}, effects of discontinuity, due to singularities \cite{j15sens2fold}, bifurcations \cite{wiercigroch10,ks10,gk10}, and now it seems also due to the hidden dynamics of the switch itself. To analyse these models and apply them to physical systems requires looking closer at how piecewise smooth models are formulated. The regularizations in \sref{sec:smooth} indicate that the range of models possible is much wider than sometimes appreciated, while the sigmoid series expansions in \sref{sec:sigseries} attempt to form a general approach to encoding these different systems via nonlinearity. 

In conclusion, whether we model a system as smooth or nonsmooth is not the crucial question, as was illustrated in \sref{sec:ueda} in particular. What matters is that hidden terms -- those which are negligible away from the switching surface -- are respected in whichever modeling framework we choose. \Sref{sec:sigseries} provides the means to preserve these hidden terms in a piecewise smooth model, and to discover them from small deviations in the `real' system if it is actually smooth. If we choose to analyze a system using piecewise smooth theory, then \sref{sec:non} outlines the basic approach, and the examples throughout this paper are merely a glimpse of the kind of strange behaviours that remain to be explored, and which may perhaps offer new means to model the dynamics of irregular systems.



\section*{Acknowledgements}

MRJ's research is supported by EPSRC Fellowship grant EP/J001317/2.


\bibliography{../../grazcat}
\bibliographystyle{plain}

\end{document}